\RequirePackage{fix-cm}

\let\oldvec\vec
\documentclass[smallextended]{svjour3}       
\let\vec\oldvec

\smartqed  
\usepackage{graphicx}
\usepackage{amsfonts} 
\usepackage{amssymb}
\usepackage[all]{xy}
\usepackage{MnSymbol}
\usepackage{tikz}
\usetikzlibrary{trees}
\usetikzlibrary{arrows,shapes,snakes,automata,backgrounds,petri}
\usepackage[normalem]{ulem}
\usepackage{mathtools}
\usepackage[T1]{fontenc}
\usepackage{mathrsfs}
\usepackage{cite}
\usepackage{appendix}
\usepackage{enumerate}
\usepackage{comment}
\usepackage{mathrsfs} 
\usepackage{url}
\usepackage{color}
\usepackage{blkarray}
\usepackage{hyperref,pgf}
\usepackage{setspace}

\newcommand{\rnet}{\{$\mathscr{S,C,R}\}\,$}
\newcommand{\RS}{\ensuremath{\mathbb{R}^{\mathscr{S}}}}
\newcommand{\PS}{\ensuremath{\mathbb{R}_+^{\mathscr{S}}}}

\newcommand{\PbarS}{\ensuremath{\overline{\mathbb{R}}_+^{\mathscr{S}}}}
\newcommand{\kinsys}{\{$\mathscr{S,C,R,K}\}$}

\newcommand{\scrS}{\ensuremath{\mathscr{S}}}
\newcommand{\scrC}{\ensuremath{\mathscr{C}}}
\newcommand{\scrR}{\ensuremath{\mathscr{R}}}

\newcommand{\scrK}{\ensuremath{\mathscr{K}}}

\newcommand{\rxn}{\ensuremath{y \to y'}}
\newcommand{\supp}{\ensuremath{\mathrm{supp} \,}}

\newcommand{\Span}{\ensuremath{\mathrm{span} \,}}

\newtheorem{rem}[theorem]{Remark}

\usepackage{amsmath}

\DeclareMathOperator{\im}{im} 

\DeclareMathOperator{\e}{e} 
\DeclareMathOperator{\diag}{diag}

\newcommand{\RR}{{\mathbb R}}

\newcommand{\mal}{\circ}
\newcommand{\st}{\mid} 
\newcommand{\dd}[2]{\frac{\text{d} #1}{\text{d} #2}}
\newcommand{\DD}[2]{\frac{\partial #1}{\partial #2}}

\newcommand\refsig{(\hspace{-0.5pt}{\it sig}) }
\newcommand\refinj{({\it inj}\hspace{1pt}) }

\newcommand\reflin{({\it lin}\hspace{1pt}) }
\newcommand\refjac{({\it jac}) }
\newcommand\refi{($i$) }
\newcommand\refii{($ii$) }

\begin{document}

\title{Some consequences of thermodynamic feasibility for chemical reaction networks}

\subtitle{Considering thermodynamic feasibility in current CRN research.}


\author{Gunter Neumann  
}


\institute{G. Neumann \at
              FAU Erlangen \\
              Tel.: +49-176-39341785\\
              \email{gunter.neumann@gmail.com}           
}

\date{Received: date / Accepted: date}

\maketitle

\begin{abstract}

Power law dynamics is used to describe the stability behavior in metabolic networks such as chemical reaction networks (CRN's). These systems allow multiple steady states within a single stoichiometric class. On the other side thermodynamic constraints such as loop-less fluxes represented by the Gorban theorem of alternatives  applied to these networks reveal considerable restrictions to their dynamics by eliminating multistability of CRN's in general. 
Thermodynamic feasible CRN's  are contained  in the class of injective CRN's. We can give an alternative proof of the  detailed balance with Brewer's Fixed Point Theorem. Furthermore we can derive by the loop-less principle the extended detailed balance. This paper establishes a link between recent research in CRN theory and thermodynamic basics.
The result has also consequences for the picture of multiple steady states as assumed for cell differentiation and regulation. CRN's provide from their perspective not enough means to maintain multistability without regulation or external control.  
\end{abstract}

\section{Introduction}
\label{intro}
Your text comes here. Separate text sections with
The dynamics of biochemical reaction networks are tightly constrained by fundamental thermodynamic laws, which impose well defined relations between reaction rate constants within that systems. Reconstructing biological reaction  networks from empirical data do have parameters that are not in accordance with these thermodynamic constraints. A systematic description of chemical reaction networks  is needed to avoid erronous and contradicting models. In \cite{gilles2007} a thermodynamic feasible model for chemical reaction networks is given. The constraints given there for chemical reaction constants does not give immediate patterns on the fluxdistribution on the whole reaction network under consideration \cite{goutsias2011}. Another approach to circumvent that high dimensional space for the choice of thermodynamical feasible reaction rates is given in 
\cite{noorlewis,palsson,beard}. There a dual approach is applied for the definition of thermodynamical feasible fluxpatterns by using the kernel of the stoichiometric matrix.  In \cite{noorlewis} thermodynamic feasible fluxes are orthogonal to the kernel of the stoichiometric matrix. 

On the other side we have an enduring  discussion about situations where chemical reacion networks are injective or multistable as introduced in  \cite{CraciunFeinberg2005}. In  \cite{mullerfeliu} conditions are given for injectivity. The loop-less and so called thermodynamically feasible fluxes outlined and specifed in \cite{noorlewis} obey in an almost natural way the injectivity conditions in \cite{mullerfeliu}. There is a long history of achievements analyzing injectivity and multistationarity of chemical reaction networks (CRN's) (\cite{Feinberg1987}, \cite{feinberg1995},\cite{HornJackson1972}, \cite{soule},\cite{thomas81}). There have been numerous refinements and generalizations of prevoius resuslts in  (\cite{CraciunFeinberg2005},\cite{mullerfeliu},\cite{angelipetrinet},\cite{banajicraciun},\cite{conradiflockerziJMB2012},\cite{joshishiu},\cite{shiusturmfels}). We would like to insert thermodynamical requirements \cite{noorlewis} into CRN's as  manifestet in \cite{mullerfeliu} to elucidate their consequences for their stability behavior. 

We will analyse thermodynamic feasible CRN's with respect to internal fixed points. 
We allow CRN's obeying thermodynamic feasibility. 
Injektivity here is defined that all participating chemical reactant species are nonvanishing at equilibrium points, i.e.  there is an internal fixed point for a CRN. An extension to the extended detailed balance \cite{Gorban2011}  will be given. The paper will be mostly self contained following other cited papers. The reader familiar with CRN notations and proofs can skip Sections \ref{sec_1} to \ref{InjectivityPol} and start with parts of Section \ref{thermodynamicBASICS}.

\section{Background material}
\label{sec_1}
In this section we give basic notations of chemical reaction networks and thermodynamic constraints. We will follow notations as used in \cite{feinberg_lectureschemical_1979,CraciunFeinberg2005,ShinarFeinberg2012,mullerfeliu,noorlewis}.

\section{Chemical reaction network theory notations}\label{FeinbergCRN}
\label{sec:CRNTPreliminaries}

Here we provide adapted reaction network theory concepts required for defining and proving the results in this article. We will generalize the reaction network theory in order to prove injectivity for classes of such networks.

\subsection{Notation}

Here we summarize the basic notations used throughout this work.
Moreover, we elaborate on the concept of sign vectors.

We denote by $I$  a finite Index set representing a collection of  $\vert I \vert =n$ species (nodes) in a network. 
Every species is determined by its  concentration $\{x_i \in \mathbb{R} \vert i \in I \}$ and we therefore denote the vector space of real-valued functions with domain $I$ by $\mathbb{R}^I \approx \mathbb{R}^n$.  
The subset of $\mathbb{R}^I$ consisting of only positive (non-negative) coordinate values is denoted $\mathbb{R}_+^I\: (\overline{\mathbb{R}}_+^I)$. For $x \in \mathbb{R}^I$ and $i \in I$,  the symbol $x_i$ denotes the value (concentration) assigned to species $i$. 
For any natural number $n \in \mathbb{N}$, we define $[n]=\{1,\dots,n\}$. 
We will use the notation $I=\{1,2,3, \cdots , n\}=[n]$ if there is no other explicit definition given.

We denote the strictly positive real numbers by $\RR_+$
and the non-negative real numbers by $\overline{\RR}_+$. For $x \in \RR^n$ we denote $x > 0$ if $x_i>0$ for all $i \in [n]$.
We define $\e^x \in \RR^n_+$ for $x \in \RR^n$ component wise, that is, $(\e^x)_i = \e^{x_i}$;
analogously, $\ln(x) \in \RR^n$ for $x \in \RR^n_+$ and $x^{-1} \in \RR^n$ for $x \in \RR^n$ with $x_i \neq 0$.
For $x,y \in \RR^n$,
we denote the component wise (or Hadamard) product by $x \mal y \in \RR^n$, that is, $(x \mal y)_i = x_i y_i$.
Further,
we define $x^b \in \RR$ for $x \in \RR^n_+$ and $b \in \RR^n$ as $x^b = \prod_{i=1}^n x_i^{b_i}$.

Given a matrix $B \in \RR^{r \times n}$, we denote by $b^1, \dots, b^n$ its column vectors
and by $b_1, \dots, b_r$ its row vectors.
Thus, the $j$th coordinate of the map $x^B \colon\RR^n_+ \to \RR^r_+$
is given by
\[
(x^B)_j  = x^{b_{j}} =x_1^{b_{j1}} \cdots x_n^{b_{jn}}.
\]
Further we define $B_\lambda$ for $B \in \RR^{r \times n}$ and $\lambda \in \RR_+^n$
as $B_\lambda = B \diag(\lambda)$.

We identify a matrix $B \in \RR^{r \times n}$ with the 
corresponding linear map $B \colon \RR^n \to \RR^r$
and write $\im(B)$ and $\ker(B)$ for the respective vector subspaces.
For a subset $S \subseteq \RR^n$,
we write $S^* = S \setminus \{0\}$
and denote the image of $S$ under $B$ by
\[
B(S) = \{ B \, x \st x \in S \}.
\]
Given sets $I \subseteq [n]$ and $J\subseteq [r]$,
we denote the submatrix of $B$ with row indices in $J$ and column indices in $I$ by $B_{J,I}$.

Whenever $U$ is a linear subspace of $\mathbb{R}^I$, we denote by $U^\perp$ the orthogonal complement of $U$ in $\mathbb{R}^n$ with respect to the standard scalar product.

By the \emph{support} of $x \in \mathbb{R}^I$, denoted $\supp(x)$, we mean the set of indices $i \in [n]$ for which $x_i \neq 0$.  When $\xi$ is a real number, the symbol $\sigma(\xi)$ denotes the sign of $\xi$.

 \begin{definition} \label{def:main_sig} 
For $x \in \mathbb{R}^n$, $\sigma(x)$ denotes the \emph{sign vector} with domain $[n]$ defined by \[\sigma(x)_i := \sigma(x_i) \in \{-1,0,1\}, \hspace{0.2 cm} \forall i \in [n]. \]

\end{definition}
Note that a sign vector $\nu\in \{-,0,+\}^n$ corresponds to the (possibly lower-dimensional) orthant of $\RR^n$ 
 given  by $\sigma^{-1}(\nu)$.
 For a subset $S \subseteq \RR^n$, we write $\sigma(S) = \{ \sigma(x) \st x \in S \}$ for
      the set of all sign vectors of $S$.

	For $x,x' \in \mathbb{R}^n$ we say that  the sign pattern of $x$ is contained in the sign pattern of $x'$ in notation \[\sigma(x) \subseteq \sigma(x') \] if we have $supp(x)\subseteq supp(x')$ where $\sigma(x_i)=\sigma(x'_i)$ $\forall \ i \in supp(x)$.

Now we are ready to state some consequences of Definition \ref{def:main_sig}.

For $x,y \in \RR^n$, we have the equivalence
\[
\sigma(x) = \sigma(y)
\quad \Leftrightarrow \quad
x = \lambda \mal y \textrm{ for some } \lambda \in \RR^n_+ ,
\]
and hence, for $S \subseteq \RR^n$, we obtain
\begin{equation}\label{eq:SigmaS}
\Sigma(S) = \sigma^{-1}(\sigma(S)) = \{\lambda \mal x \st \lambda \in \RR^n_+ \text{ and } x \in S\}.
\end{equation} 
for the union of all (possibly
      lower-dimensional) orthants that $S$ intersects. 
 For convenience, we introduce $S^* = S \setminus \{0\}$.
 
For subsets $X,Y \subseteq \RR^n$, we have the equivalences
\begin{equation} \label{sigmaSigma}
\Sigma(X) \cap Y = \emptyset
\quad \Leftrightarrow \quad
\sigma(X) \cap \sigma(Y) = \emptyset
\quad \Leftrightarrow \quad
X \cap \Sigma(Y) = \emptyset .
\end{equation}

We will now  give component wise operations for vectors $x \in \mathbb{R}^{n}$.

\noindent	

Analogously we can define a standard basis: For each $i \in [n]$, we denote by $\omega_i$ the element of $\mathbb{R}^n$ such that $(\omega_i)_j = 1$ whenever $j=i$ and $(\omega_i)_j = 0$ whenever $j \neq i$.

The \emph{standard basis} for $\mathbb{R}^n$ is the set $\left\{\omega_i \in \mathbb{R}^n : i \in [n]\right\}$. Thus, for each $x \in \mathbb{R}^n$, we have the representation $x = \sum_{i \in [n]} x_i \omega_i$. The \emph{standard scalar product} in $\mathbb{R}^n$ is defined as follows: If $x$ and $x'$ are  elements  of $\mathbb{R}^n$, then \[ x \cdot x' = \sum_{i \in [n]} x_i x'_i .\] 

 It will be understood that $\mathbb{R}^n$ carries the standard scalar product and the norm derived from the standard scalar product.  $\mathbb{R}^n$ carries the corresponding norm topology.


\subsection{Some definitions and working example}

	We will use the reaction network displayed in \eqref{EQ:ExampleNetwork} to motivate some of our definitions. Following Horn and Jackson \cite{HornJackson1972}, we call the objects at the heads and tails of the reaction arrows --- $2A, B, C, C+D$ and $E$ in \eqref{EQ:ExampleNetwork} --- the \emph{complexes} of the network. In this way, we can view the network as a directed graph, with complexes playing the role of the vertices and reaction arrows playing the role of the edges.
	
\begin{example}\label{workingexample}
\begin{align}
\label{EQ:ExampleNetwork}
\nonumber &\hspace{-2mm}2A \hspace{2 mm} \rightleftarrows \hspace{1.5 mm} B \\
\nonumber & \hspace{.5 mm} \updownarrows \hspace {3 mm}  \swarrow \\
& \hspace{1 mm} C \\[0.75 em]
\nonumber & C + D \rightleftarrows E
\end{align}
\end{example}

\begin{rem}
We will be now working in \RS, where \scrS \; is the set of species in a network which we already indicated by the equivalence $\scrS \triangleq I$ where $| \scrS | = n$. In this special case, it is advantageous to replace symbols for the standard basis of \RS  \  with the names of the species themselves. Thus, if the species in the network are given by $\scrS = \{A,B,C,D,E\}$ then a vector such as $\omega_C +\omega_D \in \RS$ can instead be written as $C+D$, and $2\omega_A$ can be written as $2A$. In fact, \RS\, can then be identified with the vector space of formal linear combinations of the species.  In this way, the complexes of a reaction network with species set \scrS \;can be identified with vectors in \RS as linear combinations of basis vectors $w_i$ for $i \in \scrS$ which we interchangeably identify with the set $[n]$ or in the case of Example \ref{workingexample} with the set of roman letters $\{A,B,C, \cdots \}$.

\end{rem}

\begin{definition}
\label{DEF:ChemicalReactionNetwork}
A \emph{chemical reaction network} consists of three finite sets:
\begin{enumerate}
	\item a set $\mathscr{S}$ of distinct \emph{species} of the network with cardinality $| \scrS | = n$;
	\item a set $\mathscr{C} \subset \overline{\mathbb{R}}_+^{\mathscr{S}}$ of distinct \emph{complexes} of the network with cardinality $| \mathscr{C} | = r$;
	\item a set $\mathscr{R} \subset \mathscr{C} \times \mathscr{C}$ of distinct \emph{reactions}, with  $| \mathscr{R} | = p$ and following properties:
		\begin{enumerate}
			
			\item $(y,y) \notin \mathscr{R}$ for any $y\in \mathscr{C}$ excluding identity reactions;
			\item for each $y \in \mathscr{C}$ there exists $y' \in \mathscr{C}$ such that $(y,y') \in \mathscr{R}$ or such that $(y',y) \in \mathscr{R}$ excluding complexes which are not participating in any of the reactions $\mathscr{R}$.
			
\item A reaction $y \to y' \in \mathscr{R}$ will be denoted  by the ordered tuple $(y,y') \in \mathscr{R}$.			
			
		\end{enumerate}
	\end{enumerate}
\end{definition}
\smallskip

If $(y,y')$ is a member of the reaction set $\mathscr{R}$, we say that $y$ \emph{reacts to} $y'$, and, following the usual notation in chemistry, we write $y \rightarrow y'$ to indicate the reaction whereby complex $y$ reacts to complex $y'$. We call the complex situated at the tail of a reaction arrow the \emph{reactant complex} of the corresponding reaction, and the complex situated at the head of a reaction arrow the reaction's a \emph{product complex}.
\smallskip

	The set of species of the network depicted in \eqref{EQ:ExampleNetwork} is $\scrS = \{A, B, C, D, E\}$. The set of complexes of the network is $\scrC = \{2A, B, C, C + D, E\}$. The set of reactions of the network is $\scrR = \{2A \to B, B \to C, C \to 2A, 2A \to C, C + D \to E, E \to C + D\}$.
\smallskip

	The diagram in \eqref{EQ:ExampleNetwork} is an example of a \emph{standard reaction diagram}: each complex in the network is displayed precisely once, and each reaction in the network is indicated by an arrow in the obvious way. Due to historical reasons we will refer to the \emph{linkage classes} of a reaction network, which can be identified with the connected components of the network whose vertices are complexes ($\mathscr{C}$) and directed edges between the complexes representing the reactions ($\mathscr{R}$). Thus, in network \eqref{EQ:ExampleNetwork} there are two linkage classes, containing, respectively, the complexes $\{2A, B, C\}$ and $\{C + D, E\}$. For a formal definition of a linkage class see \cite{feinberg_lectureschemical_1979}. We will refer to that definition  in the context of complex balancing and species balancing. 

In that context we introduce the idea of \emph{weak reversibility} which resides on the concept of linkage classes. The following definition provides some preparation.    

\begin{definition}
\label{DEF:UltimatelyReactsTo}
A complex $y \in \mathscr{C}$ \emph{ultimately reacts} to a complex $y' \in \mathscr{C}$ if any of the following conditions is satisfied:
\begin{enumerate}
\item{$y \rightarrow y' \in \mathscr{R}$};
\item{There is a sequence of complexes $y(1), y(2), \ldots , y(k)$ such that \[ y \rightarrow y(1) \rightarrow y(2) \rightarrow \ldots \rightarrow y(k) \rightarrow y'.\]}
\end{enumerate}
\end{definition}
\smallskip

In our example, the complex $2A$ ultimately reacts to the complex $C$, but the complex $C$ does not ultimately react to the complex $C + D$. 

\begin{definition}
\label{EQ:WeaklyReversible}
A reaction network \rnet is called \emph{weakly reversible} if for each $y, y' \in \mathscr{C}$, $y' \mbox{ ultimately reacts to } y$ whenever $y \mbox{ ultimately reacts to } y'$. A network is called \emph{reversible} if $y' \rightarrow y \in \mathscr{R}$ whenever $y \rightarrow y' \in \mathscr{R}$.
\end{definition}

Network \eqref{EQ:ExampleNetwork} is an example of a weakly reversible reaction network that is not reversible. Note that any  reversible  network is also weakly reversible. Note also that whenever a weakly reversible reaction network is displayed as a standard reaction diagram, every arrow in the diagram belongs to  a directed cycle of arrows.
\smallskip

\begin{definition}
\label{DEF:Reaction vectors}
The \emph{reaction vectors} for a reaction network \rnet\  are the members of the set of differences between reactant and product complex of a reaction \[ \left\{[y' - y] \in \RS : y \rightarrow y' \in \mathscr{R} \right\}. \] We will also denote   the reaction vectors of both complexes of a reaction by 
\[\Delta y = [y-y'] \] enumerated with the index set $\mathscr{R}$ :
\[\Delta y (i) \in \mathbb{R}^{\mathscr{S}} \setminus \{0\} , \ \forall \ i \in \mathscr{R} .\]
\end{definition}
\smallskip

In our example the reaction vector corresponding to the reaction $2A \to B$ is $B - 2A$, the reaction vector corresponding to the reaction $C + D \to E$ is $E - C - D$, and so on.

\begin{definition}
\label{DEF:StoichiometricSubspace}
The \emph{stoichiometric subspace} $S$ of a reaction network \rnet \  is the linear subspace of $\RS$ defined by 

\begin{equation}
S := \Span \left\{ [y' - y] \in \RS : y \rightarrow y' \in \mathscr{R} \right\}.
\label{EQ:StoichiometricSubspace}
\end{equation}
\end{definition}
\smallskip

We note that, for a reaction network \rnet, the stoichiometric subspace $S$ will often be a proper subspace of $\RS$. In other words, it will often be the case that the dimension of $S$ will be smaller than the number of species in the network. For example, in network \eqref{EQ:ExampleNetwork} we have $\dim S = 3$ while $\dim \RS = | \scrS | = 5$. In fact, the stoichiometric subspace will be a proper subspace of \RS \, whenever the network is \emph{conservative}:

\begin{definition}
\label{DEF:Conservative}
A reaction network \rnet \ is \emph{conservative} whenever the orthogonal complement $S^\perp$ of the stoichiometric subspace $S$ contains a strictly positive member of $\RS$: \[S^\perp \cap \PS \neq \emptyset.\]
\end{definition}

Network \eqref{EQ:ExampleNetwork}, for example, is conservative: it is easy to verify that the strictly positive vector $(A + 2B + 2C + D + 3E) \in \PS$ is orthogonal to each of the reaction vectors of \eqref{EQ:ExampleNetwork}.

 We will assume mass conserving reaction networks.  An intuitive example where a conservative reaction network occurs is the existence of a positive vector which we denote a (\emph{mass density vector})
\[c_m \in \mathbb{R}^{\mathscr{S}}_+\]
representing the mass $c_m(i)$ of each species $i \in \scrS$. For every reaction vector  we have 
\[\Delta y(i) \cdot c_m = 0, \ \forall \ i \in \mathscr{R}.\] As a consequence 
\[S \perp c_m \ \ or \ \  c_m \in S^{\perp}\] and the reaction network is conservative. 
\bigskip

If \rnet\  is a reaction network, then a mixture state will generally be represented by a \emph{composition} $x \in \PbarS$, where, for each $i \in \mathscr{S}$, we understand $x_i$ to be the molar concentration of species $i$. By a \emph{positive composition} we mean a strictly positive composition --- that is, a composition in $\PS$.

\begin{definition}
\label{DEF:Kinetics}
A \emph{kinetics} $\mathscr{K}$ for a reaction network \rnet\ is an assignment to each reaction $y \to y' \in \scrR$ of a \emph{rate function} $\mathscr{K}_{y \to y'}:\PbarS \to \overline{\mathbb{R}}_+$ such that \[ \mathscr{K}_{y \to y'}(x) > 0 \mbox{ if and only if } \supp y \subset \supp x.\]
\end{definition}
\smallskip

\begin{definition}
\label{DEF:KineticSystem}
A \emph{kinetic system} \kinsys\  is a reaction network \rnet\  taken with a kinetics $\mathscr{K}$ for the network. 
\end{definition}
\smallskip

\begin{definition}
\label{DEF:MassActionKinetics}
A kinetics $\scrK$ for a reaction network \rnet is \emph{mass action} if, for each reaction $y \to y' \in \scrR$, there is a positive number $\kappa_{y \to y'}$ such that the rate function $\scrK_{\rxn}$ takes the form
\begin{equation}
\label{EQ:MassActionKinetics}
\scrK_{y \to y'}(x) = \kappa_{y \to y'} x^y.
\end{equation}

\noindent The positive number $\kappa_{y \to y'}$ is the \emph{rate constant} for reaction $y \to y'$. 
We denote the vector of all rate constants by $\kappa = \{\kappa_{\rho}\}_{\rho \in \scrR} \in \mathbb{R}^{\scrR} $.
The rate function 
$\scrK_{y \to y'}(x) $ from eqn. (\ref{EQ:MassActionKinetics})   can be also identified with the \emph{flux}  $\nu_{y \to y'}$  generated by the reaction $y \to y' \in \mathscr{R}$ giving the fraction of  change of  the species concentrations (compositions) $x_i$, $i \in \mathscr{S}$ induced by the corresponding reaction vector $\Delta y$. We will distinguish between unidirectional and bidirectional fluxes.

We will also abbreviate a  \emph{kinetic system}  \kinsys\ by a chemical reaction network (CRN) if the notation is given from the context.

\end{definition}
\smallskip

\begin{definition}
\label{DEF:MassActionKineticSystem}
A \emph{mass action kinetic system} is a reaction network taken together with a mass action kinetics for the network.
\end{definition}
\smallskip

\begin{definition}
\label{DEF:SpeciesFormationRateFunction}
The \emph{species formation rate function} for a kinetic system \kinsys\  with stoichiometric subspace $S$ and rate constants $\kappa$ is the map $f_{\kappa}:\PbarS \to S$ defined by

\begin{equation}
\label{EQ:SpeciesFormationRateFunction}
f_{\kappa}(x) = \sum_{\rxn \in \scrR} \scrK_{y \to y'} (x) [y' - y]. 
\end{equation}
\end{definition}
\smallskip
For a kinetic system \kinsys\  whose underlying network is our example network \eqref{EQ:ExampleNetwork}, the species formation rate function has the following species-wise form:

\begin{align}
\label{EQ:CoordinateForm}
&f_A(x) = -2\scrK_{2A \to B}(x) +2\scrK_{B \to 2A}(x) + 2\scrK_{C \to 2A}(x) - 2\scrK_{2A \to C}(x), \\
\nonumber &f_B(x) =-\scrK_{B \to C}(x) + \scrK_{2A \to B}(x) - \scrK_{B \to 2A}(x) , \\
\nonumber &f_C(x) =  \scrK_{2A \to C}(x) -\scrK_{C \to 2A}(x) +\scrK_{B \to C}(x) -\scrK_{C + D \to E}(c) + \scrK_{E \to C + D}(c), \\
\nonumber &f_D(x) = -\scrK_{C + D \to E}(x) + \scrK_{E \to C + D}(x), \\
\nonumber &f_E(x) = -\scrK_{E \to C + D}(x) + \scrK_{C + D \to E}(x).
\end{align}

\begin{definition}
\label{DEF:DifferentialEquation}
The \emph{differential equation} for a kinetic system with species formation rate function $f_{\kappa}(\cdot)$ is given by

\begin{equation}
\label{EQ:DifferentialEquation}
\dot{x} = f_{\kappa}(x).
\end{equation} 

\end{definition}
\smallskip

Let \kinsys\  be a kinetic system. From equations \eqref{EQ:StoichiometricSubspace}, \eqref{EQ:SpeciesFormationRateFunction}, and \eqref{EQ:DifferentialEquation} we observe that the vector $\dot{x}$ will invariably lie in the stoichiometric subspace $S$ of the network \rnet. Thus, the difference of any two compositions $x \in \PbarS$ and $x' \in \PbarS$ that lie  along   the same solution of (\ref{EQ:DifferentialEquation}) will always reside in $S$. This motivates the following definition: 

\begin{definition}
\label{DEF:StoichiometricallyCompatible}
Let \rnet\  be a reaction network with stoichiometric subspace $S$. Two compositions $x$ and $x'$ in $\PbarS$ are called \emph{stoichiometrically compatible} if $x' - x\in S$.
\end{definition}
\smallskip

We note that stoichiometric compatibility is an equivalence relation. As such, it partitions $\PbarS$ into equivalence classes that we call  \emph{stoichiometric compatibility classes}.   Thus, the stoichiometric compatibility class containing an arbitrary composition $x_0$, denoted $(x_0 + S) \cap \PbarS$, is given by

\begin{equation}
\label{EQ:StoichiometricCompatibilityClass}
(x_0 + S) \cap \PbarS = \left\{ x' \in \PbarS: x' - x_0 \in S \right\}.
\end{equation}
We observe, as the notation suggests, that $(x_0 + S) \cap \PbarS$ is the intersection of $\PbarS$ with the translated hyper surface  $S$ containing $x_0$ which we denote now by $S_{x_0}=(x_0 + S) \cap \PbarS$. 

For any initial value $x_0\in \RR^n_+$
the solution is confined to coset $x_0+S_\kappa$, where $S_\kappa$
 is the smallest vector subspace containing
  the image of $f_\kappa$. We are in general interested in the positive solutions
  to the equation $f_\kappa(x) = 0$ (\ref{EQ:DifferentialEquation}) within cosets $x^\prime + S_\kappa$ with $x^\prime \in \RR^n$. 
  Due to the form of $f_\kappa$, 
  one has $S_\kappa \subseteq S$. In many applications, 
$S_\kappa = S$ for all $\kappa \in \RR_+^r$.
  If $f_\kappa$ is injective on  $(x^\prime + S)\cap \RR^n_+$ (S-injective), then 
  $f_\kappa(x) \neq f_\kappa(y)$ for all distinct $x$, $y\in (x^\prime
  +S)\cap \RR^n_+$, and, hence, the coset $x^\prime +S$ contains at most
  one positive steady state. Clearly, for a
  vector subspace $S$ of $\RR^n$, two vectors $x$, $y \in \RR^n$ lie
  in $x^\prime +S$ for some $x'\in \RR^n$, if and only if $x-y\in S$.

	A stoichiometric compatibility class will typically contain a wealth of (strictly) positive compositions. We say that a stoichiometric compatibility class is \emph{nontrivial} if it contains a member of \PS. To see that a stoichiometric compatibility class can be trivial, consider the simple reaction network \mbox{$A+B \rightleftarrows C$}, and let $\bar{x}$ be the composition defined by $\bar{x}_A = 1,  \bar{x}_B = 0, \bar{x}_C = 0$. Then the stoichiometric compatibility class containing $\bar{x}$ has $\bar{x}$ as its only member.

\begin{rem}
If a reaction network \rnet \ is mass preserving with a mass density vector $c_m \in \mathbb{R}^{\mathscr{S}}_+$ then it is conservative. We can see that $S_{x_0}$ is a compact subspace of $\overline{\mathbb{R}}^{\mathscr{S}}_+$. $S_{x_0}$ is contained in a hyper surface of the form  
\begin{equation}
S_{x_0} \subseteq c_m^{\perp}(x_0)= \{x_0 +x \in \overline{\mathbb{R}}^{\mathscr{S}}_+ | x \cdot c_m =0  \}
\end{equation} 
which is spanned between the vertices with the $i$-th coordinate set to  $\bar{x}_i=\frac{x_0 \cdot c_m}{c_m(i)}$, $\forall \ i \in \mathscr{S}$ and for all  others $j\neq i$ set to zero. 
\end{rem}

\begin{definition}
\label{DEF:Equilibrium}
\noindent An \emph{equilibrium} of a kinetic system \kinsys\  is a composition $x \in \PbarS$ for which $f(x) = 0$. A \emph{positive equilibrium} of a kinetic system \kinsys\  is an equilibrium that lies in $\PS$.
\end{definition}
\smallskip

	In light of Definition \ref{DEF:Kinetics}, a kinetic system can admit a positive equilibrium only if its reaction vectors are positively dependent:

\begin{definition}
\label{DEF:PositivelyDependent} The reaction vectors for a reaction network \rnet\  are \emph{positively dependent} if for each reaction $y \to y' \in \scrR$ there exist  positive numbers $\kappa_{y \to y'}$ such that
\begin{equation}
\sum_{\rxn \in \scrR}\kappa_{\rxn}[y' - y] = 0.
\end{equation}
\end{definition}

\begin{rem}
For any weakly reversible network, the reaction vectors are positively dependent \cite{feinberg_lectureschemical_1979}.
\end{rem}
\medskip


\section{Injectivity criterions for generalized polynomial maps}\label{InjectivityPol}

We will now follow a definition that ensure a unique fixed point in the stoichiometry class $S_{x_0}$.  For reasons of completeness and to make the paper selfcontained we follow \cite{mullerfeliu} with adaptation of notation if necessary.
We will consider families of maps defined on the positive orthant,
associated with two real matrices of coefficients and exponents, respectively, and a vector of positive parameters. We will now generalize the notation for the species formation rate function of equations (\ref{EQ:SpeciesFormationRateFunction}) and (\ref{EQ:DifferentialEquation}) of a kinetic system \kinsys\ .

\begin{definition} \label{def:main_gen}
Let $A = (a_{ij}) \in \RR^{m\times r}$, where we can identify $m=|\scrS|$  and $r= | \scrR |$. Furthermore we have $B= (b_{ij}) \in\RR^{r\times n}$, and $\kappa \in \RR_+^r$ where we can set $n \in \mathbb{N}$.
We define the associated {\em generalized polynomial map} $f_\kappa \colon \RR_+^n \to \RR^m$ as
\begin{equation*}
f_{\kappa,i}(x) = \sum_{j =1}^r a_{ij} \, \kappa_j \, x_1^{b_{j1}} \cdots x_n^{b_{jn}}, \quad i = 1, \dots, m.
\end{equation*}
\end{definition}
\begin{rem}
In that notation the species formation rate function can be identified by the equivalences $a_{ij}=\Delta
y_i(j)$ where the columns of {\emph{ stoichiometric matrix}} $A$ consists of the reaction vectors and the rows of matrix $B$ are  the reactant complexes. In our case we have the situation where $n=m$ since we assume a kinetic system  on the compositions of species $\scrS$ in a CRN.
\end{rem}
Since we had no restrictions for the coefficients of the kinetic system \kinsys\ the  term {\em generalized} indicates that we allow also polynomials with real exponents.

We often use a more compact notation.
By introducing $A_\kappa \in \RR^{m \times r}$ as $A_\kappa = A \diag(\kappa)$
and $x^B \in \RR^r_+$ via $(x^B)_j = x_1^{b_{j1}} \cdots x_n^{b_{jn}}$ for $j=1,\ldots,r$,
we can write
\begin{equation} \label{eq:map_f}
f_\kappa(x) = A_\kappa \, x^B .
\end{equation}

A generalized polynomial map $f_\kappa \colon \RR_+^n \to \RR^n$ \eqref{eq:map_f} with $A \in \RR^{n \times r}$ and $B \in \RR^{r \times n}$,
induces a system of ordinary differential equations (ODEs)
call now a {\em power-law system}:
\begin{equation} \label{eq:diffeq}
\dd{x}{t} = f_\kappa(x).
\end{equation}
We also decompose the rate of change of a composition $x$  into a multiplication of the \emph{flux} $\nu$  between the complexes $y \in \scrC$  and its stoichiometric difference matrix $A$:
\begin{equation}\label{fluxNu}
\dd{x}{t}=A_\kappa \, x^B = A \, \nu(x) \ ,
\end{equation}
where we set $\nu = \diag(\kappa) \, x^B$ for the flux. 
We subsequently set $S = \im(A)$ for the stoichiometric subspace similarly defined in  (\ref{EQ:StoichiometricSubspace}) if it is not derived from the context.
This motivates 
  the following definition of injectivity with respect to a subset.
 
\begin{definition} \label{def:main_inj}

Given two subsets $\Omega, S \subseteq \RR^n$, a function $g$ defined on $\Omega$ is called \emph{injective with respect to} $S$
if $x,y \in \Omega$, $x \neq y$, and $x-y \in S$ imply $g(x) \neq g(y)$.
\end{definition}

We will in general consider functions defined on the positive orthant, that is, $\Omega = \RR_+^n$.
When $S$ is a vector subspace, injectivity with
  respect to $S$ is equivalent to injectivity on every coset $x^\prime +S$.

Identifying $B \in \RR^{r \times n}$ with the linear map $B \colon \RR^n \to \RR^r$,
we write $B(S)$ for the image under $B$ of the subset $S \subseteq \RR^n$.
In analogy to $A_\kappa$, we introduce $B_\lambda = B \diag(\lambda)$ for  
$\lambda \in \RR^n_+$.
Finally, we write $J_{f_\kappa}$ for the Jacobian matrix associated with the map $f_\kappa$.
Here is the main result on the injectivity condition for generalized polynomial maps, which is adapted here from \cite{mullerfeliu}.

\begin{theorem} \label{thm:main}
Let $f_\kappa \colon \RR^n_+ \to \RR^m$ be the generalized polynomial map $f_\kappa(x) = A_\kappa \, x^B$,
where $A \in \RR^{m \times r}$, $B \in \RR^{r \times n}$, and $\kappa \in \RR^r_+$.
Further, let $S \subseteq \RR^n$.
The following statements are equivalent:
\begin{itemize}
\item[(inj)]
$f_\kappa$ is injective with respect to $S$, for all $\kappa \in \RR^r_+$.
\item[(jac)]
$\ker\left(J_{f_\kappa}(x) \right) \cap S^* = \emptyset$, for all $\kappa \in \RR^r_+$ and $x \in \RR^n_+$.
\item[(lin)]
$\ker(A_\kappa B_\lambda) \cap S^* = \emptyset $, for all $\kappa \in \RR^r_+$ and $\lambda\in \RR^n_+$.
\item[(sig)] $\sigma(\ker(A)) \cap \sigma(B(\Sigma(S^*))) = \emptyset$. 
\end{itemize}
\end{theorem}

\subsection{Motivation from chemical reaction networks} \label{subsec:crnt}

  For chemical reaction networks with mass-action kinetics,
  the concentration dynamics are governed by dynamical systems  \eqref{eq:diffeq} with polynomial maps $f_\kappa(x)=A_\kappa x^B$, as defined in Section \ref{sec:CRNTPreliminaries}.
We introduce some terms that are standard in the chemical engineering literature.
  One speaks of {\em
    multistationarity} if there exist a vector of rate constants $\kappa\in \RR^r_+$ and 
  two distinct positive vectors $x,y\in \RR^n_+ $ with $x-y \in S$ such that $f_\kappa(x)=f_\kappa(y)=0$. Clearly, if  $f_\kappa$ is injective with respect to $S$ for all values of $\kappa$,
  then multistationarity is ruled out.
 Therefore, Theorem~\ref{thm:main} can be applied in this setting to preclude multistationarity.

In Section~\ref{subsec:sign_monopoly},
we characterize, in terms of sign vectors,
the injectivity of a family of generalized polynomial maps with respect to a subset.
In particular,
we prove Theorem~\ref{thm:main}.
In Section~\ref{resultsCRN},
we apply our results to chemical reaction networks with power-law kinetics.
We give conditions for precluding multistationarity in the special case of a closed system as formulated in Section \ref{thermodynamicBASICS}.


\subsection{Families of linear maps} \label{subsec:sign_linear}

In this section,
we characterize, in terms of sign vectors, generalized polynomial maps $f_\kappa(x) = A_\kappa \, x^B$
that are injective with respect to a subset for all choices of the positive parameters $\kappa$.
We accomplish this through a series of results that lead to the proof of Theorem~\ref{thm:main}.


We consider the case of linear maps. We start with the following useful lemma. 

\begin{lemma} \label{lem:lin} 
Let $B \in \RR^{r \times n}$ and $S \subseteq \RR^n$. 
The following statements are equivalent:
\begin{enumerate}[(i)]
\item $\ker(B_\lambda) \cap S = \emptyset$, for all $\lambda \in \RR^n_+$.
\item $\sigma(\ker(B)) \cap \sigma(S) = \emptyset$.
\end{enumerate}
\end{lemma}

\begin{proof}
Statement \refi holds if and only if $B_\lambda \, x = B(\lambda \mal x) \not = 0$
for all $\lambda \in \RR^n_+$ and $x \in S$,
that is,
if and only if
$\ker(B) \cap \Sigma(S) = \emptyset$.
By~\eqref{sigmaSigma},
this is equivalent to statement \refii.
\end{proof}

We note that, if $0 \in S$, statements \refi and \refii do not hold, 
so we instead apply Lemma~\ref{lem:lin} to $S^*$.
In particular, if $S$ is a vector subspace of $\RR^n$,
then $\ker(B_{\lambda}) \cap S^* = \emptyset$ reduces to $\ker(B_{\lambda}) \cap S = \{0\}$,
that is, $B_{\lambda}$ is injective on $S$.

Now we  prove the equivalence of statements \reflin and \refsig in Theorem~\ref{thm:main}.

\begin{proposition} \label{pro:lin}
Let $A \in \RR^{m \times r}$, $B \in \RR^{r \times n}$, and $S \subseteq \RR^n$.
The following statements are equivalent:
\begin{enumerate}[(i)]
\item
$\ker(A_\kappa B_\lambda) \cap S = \emptyset$, for all $\kappa \in \RR^r_+$ and $\lambda\in \RR^n_+$. 
\item
$\sigma(\ker(A)) \cap \sigma(B(\Sigma(S))) = \emptyset$.
\end{enumerate}
\end{proposition}

\begin{proof}
By Lemma~\ref{lem:lin} applied to the matrix $A$ and the subset $B(\Sigma(S))$,
condition \refii is equivalent to $\ker(A_\kappa) \cap B(\Sigma(S)) = \emptyset$, 
for all $\kappa \in \RR^r_+$.
By~\eqref{eq:SigmaS}, this statement is equivalent to the fact that $\ker(A_\kappa) \cap  B_\lambda(S) = \emptyset$, 
for all  $\kappa \in \RR^r_+$ and ${\lambda \in \RR^n_+}$, which is in turn clearly equivalent to condition \refi.
\end{proof}

Again, if $S$ is a vector subspace, $\ker(A_\kappa B_\lambda) \cap S^* = \emptyset$
reduces to $\ker(A_\kappa B_\lambda) \cap S = \{0\}$,
that is, $A_\kappa B_\lambda$ is injective on $S$.
Clearly, the statements in Lemma~\ref{lem:lin} are necessary conditions for the statements in Proposition~\ref{pro:lin}.

\subsection{Families of generalized monomial/polynomial maps} \label{subsec:sign_monopoly}

In this subsection, we use   results on families of linear maps
to give sign conditions for the injectivity of families of generalized polynomial maps with respect to a subset.

We will specify Definition~\ref{def:main_inj} for our purposes here and
we conclude that a function $g$ defined on $\RR^n_+$ is injective with respect to a subset $S \subseteq \RR^n$
if and only if for every $x \in \RR^n_+$ one has $g(x) \neq g(y)$ for all $y \in S^*_{x}$.
If $S$ is a vector subspace, then $g$ is injective on the intersection $(x+S) \cap \RR^n_+$ of any coset $x+S$ with the domain $\RR^n_+$.

We continue with a key observation.

\begin{lemma} \label{lem:lambda}
For a vector subspace $S \subseteq \RR^n$, let
\begin{equation} \label{eq:Lambda}
\Lambda(S) := \{ \ln x - \ln y \st x,y \in \RR^n_+ \textrm{ and } x - y \in S \} .
\end{equation} 
Then, $\Lambda(S) = \Sigma(S)$.
\end{lemma}

\begin{proof}
Let $x,y \in \RR^n_+$ such that $x-y \in S$. Then, using the strict monotonicity of the logarithm 
we have $\sigma(\ln x - \ln y) = \sigma(x-y) \in \sigma(S)$ and  hence $\ln x - \ln y \in \Sigma(S)$.
This proves the inclusion $\Lambda(S) \subseteq \Sigma(S)$.
Conversely,
let $\lambda \in \RR^n_+$ and $z \in S$.
We construct $x,y \in \RR^n_+$ such that $\ln x - \ln y = \lambda \mal z$ 
and $x - y = z$ as follows:
if $z_i \neq 0$, then $e^{\lambda_i z_i} \neq 1$, so we may define 
$y_i := z_i / (\e^{\lambda_i z_i} - 1)$ and $x_i := y_i \e^{\lambda_i z_i}$;
otherwise, set $x_i=y_i=1$.
This proves $\Sigma(S) \subseteq \Lambda(S)$. 
\end{proof}


\begin{lemma} \label{lem:SB}
For $B \in \RR^{r \times n}$ and $S\subseteq \RR^n$, let
\begin{equation} \label{eq:SB}
S_B := \{ x^B - y^B \st x,y \in \RR^n_+ \textrm{ and } x - y \in S^* \} .
\end{equation}
Then, $\sigma(S_B) = \sigma(B(\Sigma(S^*))).$
\end{lemma}
\begin{proof}
For $x,y \in \RR^n_+$, we have $\sigma(x^B - y^B) = \sigma(B(\ln x - \ln y))$ by the strict monotonicity of the logarithm,
and hence
\begin{align*}
\sigma(S_B)
= \sigma\big(\{ B(\ln x - \ln y) \st x,y \in \RR^n_+ \textrm{ and } x - y \in S^* \}\big)
= \sigma(B(\Lambda(S^*))),
\end{align*}
using~\eqref{eq:Lambda}.
By Lemma~\ref{lem:lambda}, $\sigma(S_B) = \sigma(B(\Sigma(S^*)))$.
\end{proof}

We  can show now an important proposition for later purposes.
\begin{proposition} \label{pro:mono}
Let $B \in \RR^{r \times n}$ and $S \subseteq \RR^n$.  
Further, let $\varphi_B \colon \RR^n_+ \to \RR^r_+$ be the generalized monomial map $\varphi_B(x) = x^B$.
The following statements are equivalent:
\begin{enumerate}[(i)]
\item $\varphi_B$ is injective with respect to $S$.
\item $\sigma(\ker(B)) \cap \sigma(S^*) = \emptyset$.
\end{enumerate}
\end{proposition}

\begin{proof}
By \eqref{eq:SB}, statement \refi is equivalent to $0 \notin S_B$.
By Lemma~\ref{lem:SB}, this is in turn equivalent to $0 \notin B(\Sigma(S^*))$,
that is, $\ker(B)\cap \Sigma(S^*)=\emptyset$.
By \eqref{sigmaSigma}, this is equivalent to statement \refii.
\end{proof}

Next we prove the equivalence of statements \refinj and \refsig in Theorem~\ref{thm:main}.

\begin{proposition} \label{pro:poly}
Let $f_\kappa \colon \RR^n_+ \to \RR^m$ be the generalized polynomial map $f_\kappa(x) = A_\kappa \, x^B$,
where $A \in \RR^{m \times r}$, $B \in \RR^{r \times n}$, and $\kappa \in \RR^r_+$.
Further, let $S \subseteq \RR^n$.
The following statements are equivalent:
\begin{itemize}
\item[(inj)]
$f_{\kappa}$ is injective with respect to $S$, for all $\kappa \in \RR^r_+$.
\item[(sig)]
$\sigma(\ker (A)) \cap \sigma(B(\Sigma(S^*))) =\emptyset$.
\end{itemize}
\end{proposition}

\begin{proof}
Statement \refinj  asserts that for $x,y \in \RR^{n}_+$ with $x-y \in S^*$, 
we have $A_\kappa \, (x^B-y^B) \neq 0$ for all $\kappa \in \RR^r_+$.
This is equivalent to asserting that 
$\ker(A_\kappa) \cap S_B = \emptyset$ for all $\kappa \in \RR^r_+$, with $S_B$ as in~\eqref{eq:SB}.
By applying Lemma~\ref{lem:lin} to the matrix $A$ and the subset $S_B$,
this is in turn equivalent to $\sigma(\ker(A)) \cap \sigma(S_B) = \emptyset$.
By Lemma~\ref{lem:SB}, $\sigma(S_B) = \sigma(B(\Sigma(S^*)))$, and the equivalence to statement \refsig is proven.
\end{proof}

To prove the equivalence of statements \reflin and \refjac in Theorem~\ref{thm:main},
we will use the following observation. 

\begin{lemma} \label{lem:jac}
Let $A=(a_{ij})\in\RR^{m \times r}$, $B=(b_{ij}) \in \RR^{r \times n}$,
 $\kappa \in \RR^r_+$, 
and
$\lambda \in \RR^n_+$.
Further,
let $f_\kappa \colon \RR^n_+ \to \RR^m$ be the generalized polynomial map $f_\kappa(x) = A_\kappa \, x^B$.
Then, the sets of all Jacobian matrices $J_{f_\kappa}(x)$ and all matrices $A_\kappa B_\lambda$ coincide:
\[
\left\{J_{f_\kappa}(x)  \st \kappa \in \RR^r_+ \textrm{ and } x \in \RR^n_+ \right\}
=
\left\{ A_\kappa B_\lambda  \st \kappa \in \RR^r_+ \textrm{ and } \lambda \in \RR^n_+\right\} .
\]
\end{lemma}

\begin{proof} 
As $f_{\kappa,i}(x) = \sum_{j=1}^r a_{ij} \, \kappa_j \, x^{b_{j}}$,
the $(i,\ell)$th entry of the Jacobian matrix of $f_{\kappa}$ amounts to
\[
J_{f_\kappa}(x) _{i,\ell} = \DD{f_{\kappa,i}(x)}{x_\ell} =
\sum_{j=1}^r a_{ij} \, \kappa_j \, x^{b_{j}} \, b_{j\ell} \, x_\ell^{-1} .
\]
That is, 
$$J_{f_\kappa}(x)= A \diag(\kappa \mal x^B) B \diag(x^{-1}) = A_{\kappa'} B_\lambda$$
with $\kappa' = \kappa \mal x^B$ and $\lambda = x^{-1}$.  
Clearly, quantifying over all $\kappa \in \RR_+^r$ and $x \in \RR^n_+$ is equivalent to quantifying 
over all $\kappa' \in \RR^r_+$ and $\lambda \in \RR_+^n$.
\end{proof}

We can now combine all the results in this section in the proof of the preliminary main theorem.

\begin{proof}[Proof of Theorem~\ref{thm:main}]
The equivalences \reflin $\Leftrightarrow$ \refsig and \refinj $\Leftrightarrow$ \refsig
are shown in Propositions~\ref{pro:lin} and \ref{pro:poly}, respectively.
The equivalence \refjac $\Leftrightarrow$ \reflin follows from Lemma~\ref{lem:jac}.
\end{proof}

The case $S=\im(A)$ and $m=n$ arises in applications to chemical reaction networks,
which we need for considering explicit examples which are conservative and mass preserving.

\section{Thermodynamic basics}\label{thermodynamicBASICS}
In order to consider thermodynamic aspects in a flux distribution of  a kinetic system \kinsys\ we have to assign potential differences $\Delta{G}$ between the  complexes of each reaction of the CRN in form of a vector of  potentials for the complexes.
The Gibbs potential e.g. for reaction $A+B \leftrightarrow C+D$ is given according to equation (1)   in \cite{kuemmel} by 
\begin{equation}\label{gibbsenthalpie}
\Delta{G} = y_C G_C^0 + y_D G_D^0  - y_A G_A^0 - y_B G_B^0  + RT \, \ln(K_a)
\end{equation}
over the constant $R=N_A \cdot k_b$, the  activities 
$$K_a= x^{[y_2-y_1]}$$ 
and the zero point energies $G_x^0$ (see also \cite{kuemmel} eqn. (1)) where we assume component wise $y_1=(y_A,y_B,0,0)$ and $y_2=(0,0,y_C,y_D)$. 
By carrying out the logarithm in eqn. (\ref{gibbsenthalpie})  we can find a vector $\gamma\in \mathbb{R}^n$ for the potentials of the individual species depending on their concentrations and stoichiometric coefficient, such that we obtain 
\begin{equation}\label{gibbsdelta}
\Delta G = \gamma^T A
\end{equation}
as the differential energy between the complexes for the current temperature and species concentrations.

\subsection*{Characterization of loop-less fluxes for CRN}

The following classification of fluxes can be traced back to the Gordan theorem of alternatives \cite{noorlewis,martino2013} which we will state here: 
\begin{theorem}
(Gordan's theorem) $\forall A\in \mathbb{R}^{ n\times m}$ exactly one of the following two statements is true:
\begin{description}
\item[(a)] $\exists z \in \mathbb{R}^m_{+} \setminus \{0\}$, s.t. $Az=0$
\item[(b)] $\exists y \in  \mathbb{R}^n$ s.t. $ A^{\top} y>0$
\end{description}
\end{theorem}

In \cite{noorlewis} a transformation of the Gordan theorem for the case of reversible fluxes of a chemical reaction network for a closed system is given. 
A reaction system fully reversible will be called loop-free/thermodynamically feasible  (b) or thermodynamically not feasible with loops (a)  if the following holds:

 \begin{corollary} \label{noorgordanlemma}
 For all $\hat{A} \in \mathbb{R}^{n \times r}$ where $n$ is the number of species and $r$ the number of (bidirectional/reversible) reactions and every $\nu \in \mathbb{R}^r$ exactly one of the following cases is true:
 \begin{description}
 \item[(a)] $\exists \hat{z}  \in \mathbb{R}^r \setminus  \{0\}$, s.t. $(\forall i \  sign(\hat{z}_i) \in \{sign(\nu_i),0\})\wedge \hat{A} \hat{z}=0$
 \item[(b)] $\exists \gamma \in \mathbb{R}^n$ s.t. $(\forall i \ sign(\hat{A}^{\top}\gamma)_i = - sign(\nu_i) \vee  \nu_i=0)$
 \end{description}
\begin{proof}
See \cite{noorlewis}.
\end{proof}
\end{corollary}

The existence of a potential for the complexes in a reaction network out rules the existence of  fluxes within a closed network without sources. Therefore that alternative expresses the fact that we cannot have a flux keeping the concentrations of the species constant when there are differences between the potential of the complexes. The net energy consumption would be zero and the turnover would be non-zero which would be impossible due to the conservation of energy. The distribution of the associated potentials over the CRN does not allow thermodynamically infeasible fluxes. 

In  Corollary \ref{noorgordanlemma} we were choosing $\gamma$ instead of $y$ in order to avoid an overlap with the stoichiometry vector $y_i$ and also to give the link to the chemical potential introduced in equations (\ref{gibbsenthalpie}) and (\ref{gibbsdelta}) since $\gamma^{\top}A$ is equivalent to $A^{\top}\gamma$. {\bf (b)} in Corollary \ref{noorgordanlemma} reflects the fact that the flux $\nu_i$ is in opposite direction to  the increasing potential   between the  complexes which is given by $(A^{\top}\gamma)_i$.

We can transform that relation to our reversible system.   We set $m=2r$ the number of all unidirectional reactions in a fully reversible CRN   and order the signs of the flux $\nu \in \mathbb{R}^r$ with $sign(\nu_i)=d_i$ for $i \in [r]$ according to the first $r$ forward and $r$ backward fluxes or each reversible reaction where we have $d_i=-d_{i+r}$ and the total flux results as the sum of the forward and backward flux: $\nu_i = z_i-z_{i+r}$ for $z \in \mathbb{R}^m_{+}$. We can set up  the following result which is an equivalent formulation of loop-free fluxes  from Corollary \ref{noorgordanlemma}  for unidirectional fully reversible CRN's.
 
 \begin{corollary}\label{noorunidirectional}
  For all $A \in \mathbb{R}^{n \times m}$ where $n$ is the number of species and $m=2r$ the number of reactions and every $\nu \in \mathbb{R}^r$ exactly one of the following cases is true:
 \begin{description}
 \item[(a)] $\exists z  \in \mathbb{R}^m_{+} \setminus  \{0\}$ $\wedge$  $(\exists j \in [r]$ with $ z_j \neq z_{j+r} )$, s.t. $(\forall i \in [r] \  sign(z_i-z_{i+r}) \in \{sign(\nu_i),0\})\wedge A z=0$
 \item[(b)] $\exists \gamma \in \mathbb{R}^n$ s.t. $(\forall i (\ sign(A^{\top}\gamma)_i = -  \ sign(A^{\top}\gamma)_{i+r} =  - sign(\nu_i) ) \vee  \nu_i=0)$
 \end{description}
 \begin{proof}
 Equivalence between Corollary \ref{noorgordanlemma} and \ref{noorunidirectional} concerning {\bf (a)} can be seen by doubling the matrix $\hat{A}$ for the bidirectional case by setting $A=(\hat{A},-\hat{A})$ and also doubling the vector $\hat{z}$ by setting $z_i=\max{(\hat{z}_i,0)}$ and $z_{i+r}=-\min{(\hat{z}_i,0)}$ for $i \in [r]$.  The reverse can be done by halving $A$ to form $\hat{A}$ and by taking differences $\hat{z}_i=z_i-z_{i+r}$ for $i \in [r]$. {\bf (b)} is equivalent in both Corollaries.
\end{proof}
\end{corollary}

\begin{remark}\label{remarknonreversible}
Corollary \ref{noorunidirectional}  can be extended to the case where  reaction $\mathcal{R}_i$, $i \in [r]$ are not reversible by choosing $\nu \in \mathbb{R}^r$ such that the sign of $\nu_i$ is in accordance with the direction of the reaction $\mathcal{R}_i$.
\end{remark}

\begin{remark}\label{remarkortho}
The exclusion of the case {\bf (a)}, where $A$ describes  a closed CRN,  implies  that there is no component $x$ of $\nu$ that is in the  null space of $A$. The process of elimination of components $\nu$ in  $\ker(A)$  is a minimization of the total flux  which results that $\nu$ has to be orthogonal to the null space of $A$ for a closed system: 
\begin{equation}\label{AperpKer}
\nu \perp \ker(A) \ 
\end{equation}
which implies therefore 
\begin{equation}\label{AperpKerGenPol}
\kappa \circ x^B \perp \ker(A) \ 
\end{equation}

\end{remark}

\subsection*{Preliminary summary}
Up to that point we stated all necessary prerequisites to derive the consequences in the following results section. On the one side we have the thermodynamic (loop-less) constraints and on the other side we have the conditions for injectivity of CRN's. We will show with both input relations that a wide range of CRN's are injective. This implies that these networks do not have multiple steady states.

\section{Results}\label{resultsCRN}
\subsection{Main Theorem}
In that section we combine the established results about injectivity of polynomial maps in Theorem \ref{thm:main}  and the orthogonality  relation of Remark \ref{remarkortho}. We will use the notation of  CRN's as introduced in Section \ref{InjectivityPol}. By decomposing the stoichiometry matrix into
\begin{equation}\label{Adecomposition}
A=KE
\end{equation}
we have similar to eqn. (\ref{eq:map_f}) for the kinetic system of  the specific CRN 
\begin{equation}\label{crnfeinberg}
\frac{ {\sf d} x}{ {\sf d} t} =f_{\kappa}(x) = KE \ \kappa \circ x^B  \ .
\end{equation}
The columns of $K$  are the stoichiometry vectors  of all  $p$ complexes $ y_{j} \in \mathcal{C}$, $j \in [p]$ involved in the $r$ reactions in $\mathcal{R}$. $E$ is the incidence matrix between the interacting complexes forming the matrix $A$, which consists of all stoichiometric differences of the reacting complexes $[y_i -y_i^{'} ]$ with  $ y_i \rightarrow y_i^{'} \in \mathcal{R}$, $i \in [r]$. The  rows of $B$ of Section \ref{InjectivityPol} are all reactant complexes  of each reaction  and is therefore different to $K$ in general.  
We  have $S  \subset {\textit{im}}(A)$. 

\begin{lemma}\label{sigperp}
Let $V,W \subset \mathbb{R}^n$ be two subspaces for which $v \in V$ and $w \in W$ implies $v \bot w$ then $\sigma(V) \cap \sigma(W^*)=\emptyset $. (The converse does not hold). 
\begin{proof}
Assume there exists $v \in V$ and $w \in W$ s.t. $\sigma(v)=\sigma(w)\neq 0$ then $v \cdot w > 0$ which contradicts $v \perp w$. 
\end{proof}
\end{lemma}

We now use the relation in eqn. (\ref{AperpKer})  to show the following:

\begin{lemma}\label{orthokappaker}
For $x, \kappa \in {\mathbb{R}^n_{+}}$ and $A$, $B$ as in Proposition \ref{pro:lin}
\begin{equation}\label{feasibleorthogonal} 
\ker(A) \bot \ \kappa \circ x^B \Rightarrow \sigma(\ker(A))  \cap \sigma( x^B) =\emptyset\ ,
\end{equation} 
\begin{proof} 
From Lemma \ref{sigperp} it can be seen first that 

\begin{equation}\label{feasibleorthogonalBB} 
\ker(A) \bot \ \kappa \circ x^B \Rightarrow \sigma(\ker(A))  \cap \sigma(\kappa \circ x^B) =\emptyset\ ,
\end{equation} 
since $\kappa \in {\mathbb{R}^n_{+}}$ preserves the signature of $ x^B$ we have eq.  (\ref{feasibleorthogonal}). 

\end{proof}
\end{lemma}

With this result we can now  state our main theorem with Proposition \ref{pro:lin}:

\begin{theorem}\label{maintheoremfeasible}
For a thermodynamic feasible kinetic  system \kinsys\ as in equation (\ref{eq:diffeq}) where $n$ is the number of species with species compositions $x \in {\mathbb{R}^n_{+}}$ involved in $r$ reactions $\{\mathcal{R}_i\}_{i \in[r]}$  and  stoichiometric difference matrix $A \in \mathbb{R}^{n\times r}$ and reactant complex matrix $B \in \mathbb{R}^{r \times n}$ with reaction rates $\kappa \in \mathbb{R}^r_{+}$  and corresponding generalized polynomial map $f_{\kappa} (x): \mathbb{R}^{n}_{+} \rightarrow \mathbb{R}^n$ with $A_{\kappa}= A \ diag (\kappa)$ and $S \subseteq \RR^n$ holds
\begin{equation}\label{AperpxB}
 \sigma(B(\Sigma(S^*))) \cap \sigma(\ker(A)) = \varnothing , 
\end{equation}
 (c.f. eqn.  (\ref{AperpKerGenPol})) holds for all $x  \in \mathbb{R}^n_{+}$. Thermodynamically feasible systems are always injective.

\end{theorem}

\begin{proof} 

 From condition (\ref{AperpKer}) we also have $(\nu_1-\nu_2) \perp \ker(A)$ for thermodynamically feasible flows $\nu_1 \neq \nu_2$, which implies  $\kappa  (x^B - y^B) \perp ker(A)$.  With the notation 
  \[S_B := \{ x^B - y^B \st x,y \in \RR^n_+ \textrm{ and } x - y \in S^* \} \] 
  if follows $\kappa (S_B) \perp \ker(A)$ 
 and with  Lemmata \ref{sigperp} and \ref{orthokappaker}  we have  $\sigma(S_B) \cap \sigma(ker(A)) = \varnothing$
 Following Proposition \ref{lem:SB}, we obtain Equation (\ref{AperpxB}).  The relation holds for all $\kappa \in \mathbb{R}^r_{+}$.
 which relaxes the condition of eqn. (\ref{AperpKerGenPol}) to its sign condition.

\end{proof}


\begin{remark}
Condition (\ref{AperpxB}) is also a requirement for the underlying network as given by the stoichiometric difference matrix $A \in \mathbb{R}^{n\times r}$, which need to be thermodynamically feasible in that sense. 
\end{remark}

\begin{remark}
Relation (\ref{feasibleorthogonal}) holds for all $x \in \mathbb{R}^n_{+}$. Condition (\ref{feasibleorthogonal}) is a property of the underlying CRN and is independent of its reaction constants and strictly positive species concentrations. In theorem \ref{maintheoremfeasible} we also allow $\kappa$  for which thermodynamic feasibility might not be realistic. But we obtain in that case that thermodynamic feasible reaction systems from theorem \ref{maintheoremfeasible} are contained in the set of injective systems as characterized in theorem \ref{thm:main}. 
\end{remark}

We will now give some consequences of that result. 

\subsection{Detailed balance}

A preliminary Lemma:
\begin{lemma}\label{injectcomplexes}
For $[y_i-y_i^{'}]\in \mathcal{R}$ with $y_i \neq y_i^{'}$ at least one of the following two cases is true: 

\begin{description}\label{ydiffmapy}
\item[a)]$ y_i \cdot [y_i-y_i{'}] \neq 0$
\item[b)]$ y_i' \cdot [y_i-y_i{'}] \neq 0$
\end{description}
\begin{proof}
Assume that both are zero then we would have $0 < [y_i-y_i{'}]\cdot [y_i-y_i{'}]=0$.   
\end{proof}
\end{lemma}

\begin{corollary}[Detailed Balance]\label{detailed_balance}
For a  kinetic system of $r$ reversible reactions  with thermodynamic feasible fluxes the corresponding generalized polynomial $f_{\kappa}$ is injective and has a unique fixed point.  Further for a conservative kinetic system there exists an interior fixed point $x_0 \in \RR^n_+ \cap S_{x_0}$ s.t. $f_{\kappa}(x_0)=0$.
\begin{proof} 
We will give first an elementary proof and then an algebraic one. 
In the case of a fully reversible network, we have $S= im (A)$, because reversible reactions allow to define a vector space over the column space of $A$. 
From Lemma \ref{injectcomplexes} we can see that matrix multiplication between complexes and reaction differences do not vanish. Further more the row space of $B$ is the same as the columns space of $A$. 
We can check that by selecting a subset of reaction differences $[y_{k_i}-y_{k_i}'] \in \mathcal{R}$ for $ i \in [k]$ where $k = dim(S)=dim(im(A))$. In the same way we can select a subset of maximum $k \leq k'\leq 2k$   row vectors $\{y_{i_{k'}}''\}_{i \in [k']}$ of $B$ out of the $\{y_{k_i},y_{k_i}'\}_{i \in [k]}$  pairs  for which  $span {(\{[y_{k_i}-y_{k_i}']\}_{i \in [k]})} \subseteq span{(\{y_{i_k}''\}_{i \in [k']})}$ holds since the column space of $KE=A$ is contained in the row space of $B$.
Together with lemma \ref{ydiffmapy} we see that $S^* $ is mapped injectively into $\textit{im}(B)$, The image of a vector $\Delta y^{\star} \in S^*$ will be mapped under $AB$ to the non-scaled projection 
\[\sum_{i \in [r]} \Delta y_{k_i}^+ (\Delta y_{k_i}^+ \cdot \Delta y^{\star}) \neq 0, \ \Delta y_{k_i}^+ = [y_{k_i}-y_{k_i}']\]
where the sum is over all reversible reactions, with only one representative of the forward (+) and backward (-) reaction. 
$\Delta y^{\star} $ maps also injectively into  $im (A)=S$. Hence, condition (lin) from Theorem \ref{thm:main} is satisfied and $f_{\kappa}$ is injective. The parameters $\kappa$ and $\lambda$ do only contribute scalings  to the projection and do not alter the result. 

Another proof for the injectivity can be derived from Corollary \ref{noorunidirectional} where a Gibbs potential is possible due to the bidirectionality implicitly assumed in Eq. \ref{gibbsenthalpie}, which implies the existence of a vector $\gamma$ with $(\gamma A) \circ (\nu) > 0$. 

If the kinetic system is conservative $S_{x^{\star}} \cap \RR^n_+$ is isomorph to some simplex $\Delta_{S_x}$ in the positive orthant. According to Brewers Fixed Point Theorem \cite{Brewer1912}, there is at least one fixed point and due to injectivity maximal one in $\RR^n_+$.  We now show that fixed points at the border are not stable. We try to show that the interior fixed point is stable and the border is repulsive. We assume one species is extinct. The border face of  $\RR^n_+$ ($x_i = 0$). We consider reactions where $x_i$ is involved. Due to magnitude only one of the reversible reactions is active, where $(\Delta y_{k})_i >0$. Hence, $f_{\kappa,i}(x)>0$. We can exclude species which are not altered in any reaction. In general we can assume that if a species composition $x$ is close to the boundary of the simplex $\partial \Delta_{S_x}$ at least one of the species has a significant amount of the conservative mass and is involved in at least one of the reactions which have almost extinct species concentration. This reaction increases the amount of at least one minor represented species. By induction another low represented species is also involved in a reaction with  a significantly represented species and has, hence, a positive production rate. We can conclude that the boundary of a reversible system is repulsive. Therefore a reversible chemical reaction system has one unique interior fixed point due to injectivity. 
\end{proof}
\end{corollary}

\begin{remark}\label{Wegscheiders}
For the proof of the detailed balance we did not need the Wegscheiders conditions. An explicit calculation of the solution of a detailed balance for all reversible reactions is unique, when Wegscheiders conditions are added. (Calculation not shown here.) On the other side the Wegscheiders conditions guarantee that there exist thermodynamically feasible fluxes (in the reversible network), which is not explicitly given by the network topology. 
\end{remark}

\subsection{Gorban Yablonsky Theorem of extended detailed balance.}

We can now derive another proof for the Gorban Yablonsky Theorem \cite{Gorban2011}, which was independently derived in 2011 from the proof here. If thermodynamic feasibility is satisfied, we have a unique fixed point. 
\begin{lemma}\label{GoorbanYablonski}
If there is any non-reversible reaction $y_i \rightarrow y_i' \in \scrR$  in a closed kinetic system \kinsys\ satisfying thermodynamic feasibility then there is no interior fixed point  $x_0 \in \mathbb{R}^{\scrS}_+\cap S_{x_0}$ for the kinetics \scrK.
\end{lemma}
\begin{proof}
For an equilibrium point $x_0 \in \mathbb{R}^{\scrS}_+\cap S_{x_0}$ s.t. $f_\kappa (x_0)=0$ we would have for the flux $\nu= \diag(\kappa) x^B \neq 0$ since we have for the flux component of  the irreversible reaction $i$: $\nu_i = \kappa_i x^B_i>0$. This implies together with the steady state condition that $A \nu =0$ that there is a non-vanishing internal steady state flux $\nu$, which is a contradiction to Remark \ref{remarkortho}. 
\end{proof}

\begin{remark}
As a result from  Lemma \ref{GoorbanYablonski} there is no irreversible reaction within the face of the positive orthant  $\mathbb{R}^{\scrS}_+$ where $x_0$ is located. All irreversible reaction vectors point inside or traverse of this face. The remaining system must consist of a fully reversible network, where some species are annihilated. In this face an interior fixed point is again possible. 
\end{remark}


\section{Conclusion}

Including thermodynamic principles into CRN's as developed in \cite{noorlewis,palsson,beard} leads to a restriction of the available parameter space. Thermodynamic feasible reaction dynamics requires injective generalized polynomial maps for the dynamics of the species concentrations.  Reversible  CRN's imply injectivity. 

Loop-circuits in interacting networks are the basis of multistability \cite{soule,gaskins2009,angeli2004,mincheva2007,muras2012}. The exclusion of loops precludes multistability in a wide class of chemical reaction networks. 

The parameter set in the power law description of CRN's is tested to satisfy the conditions of equation (\ref{AperpKer}). We did not include explicit dependence of the reaction parameters $\kappa$ into the analysis.

For multistability as proposed as a key factor for cell differentiation (\cite{delbrueck1949,huang2005}) we want to conclude that metabolic networks are regulated by signal transduction and not by triggering intrinsic multistability.  Therefore we can assume or predict that mutistability is governed by regulatory mechanisms, which are not primarily subjected to stoichiometry,  power law kinetics and thermodynamic energy potentials. We suggest that multiple steady states as observed in differentiated cells are created by tuning reaction constants and trafficking of metabolites by channels. It is also possible to think about timescales of reaction constants that separates between reaction equilibrium that are fast and slow such that slow reactions provide constant rates and fast reactions reach equilibrium almost immediately.

\section*{Author's contributions}
    GN is responsible for the idea, proofs and design of the present work. All authors read and approved the final manuscript. The study was based upon published material and results as given in the references.


\begin{acknowledgements}

The work was done during our stay at the Friedrich Alexander University in Erlangen at the  Department of Mathematics. We would like to thank Gerhard Keller and Andreas Knauf for discussions including corrections,  critical comments and encouragements.  The work was finished after a long break in 2020. 

\end{acknowledgements}

%
\section*{Conflict of interest}
 The authors declare that they have no conflict of interest.



\end{document}